\documentclass[10pt]{amsart}
\usepackage{amssymb}
\numberwithin{equation}{section}

\newtheorem{theorem}{Theorem}[section]

\newtheorem{proposition}[theorem]{Proposition}
\newtheorem{lemma}[theorem]{Lemma}
\theoremstyle{definition}
\newtheorem{definition}[theorem]{Definition}
\theoremstyle{remark}
\newtheorem{remark}[theorem]{Remark}

\newtheorem{question}{Question}
\newcommand{\R}{\mathbb{R}}
\newcommand{\C}{\mathbb{C}}
\newcommand{\Z}{\mathbb{Z}}

\newcommand{\bbT}{\mathbb{T}}

\newcommand\lie[1]{\mathfrak{#1}}

\newcommand{\fh}{\lie{h}}
\newcommand{\fg}{\lie{g}}

\def	\inv	{^{-1}}

\newcommand{\GL} {\text{GL}}

\begin{document}

\title{Maximal tori in the contactomorphism groups of circle bundles over 
Hirzebruch surfaces}

\author{Eugene Lerman}
\address{Department of
Mathematics, University of Illinois, Urbana, IL 61801}
\email{lerman@math.uiuc.edu}

\thanks{Supported by  NSF grant DMS-980305 and by R. Kantorovitz.}
%\date{\today}

\begin{abstract}
In a recent preprint Yael Karshon showed that there exist
non-conjugate tori in a group of symplectomorphisms of a Hirzebruch
surface.  She counted them in terms of the cohomology class of the
symplectic structure.  We show that a similar phenomenon exists in the
contactomorphism groups of  pre-quantum circle bundles over 
Hirzebruch surfaces.  Note that the contact structures in question are
fillable.  This may be contrasted with an earlier paper where we
showed that there are infinitely many non-conjugate tori in the
contactomorphism groups of certain overtwisted lens spaces (Contact
cuts, {\em Israel J.\ Math}, {\bf 124} (2001), 77--92).
\end{abstract}

\maketitle

\section{Introduction}
Consider the product of two projective spaces $M=\C P^1 \times \C
P^1$. Let $\omega_1$ and $\omega_2$ denote the pull-backs of the
standard area form on $\C P^1 = S^2$ by the two projection maps (so
that $\int _{\C P^1 \times \{ *\}} \omega _1 = 1$ etc.).  Let $\omega
_{a,b} = a\omega_1 + b \omega_2$; it is a symplectic form if $a,b >
0$.  Gromov in his fundamental paper \cite{Gr} showed that the
topology of the group of symplectomorphisms of $(M, \omega_{a,b})$
changes as the ratio $a/b$ crosses integers.  The rational cohomology
ring of the connected component of the group of symplectomorphisms of
$(M, \omega_{a,b})$ has been computed for different values of $a/b$ by
Abreu \cite{Ab} and Abreu - McDuff \cite{AbMc}.  Karshon \cite{K}
showed that for $a\geq b > 0$ the number of conjugacy classes of
maximal tori in the group of Hamiltonian symplectomorphism of $(M,
\omega_{a, b})$ is the number of integers $k$ satisfying $0 \leq k <
a/b$.  Similar results hold for a one point blow-up $\widetilde{\C
P}{} ^2$ of $\C P^2$.

Now suppose that $a$ and $b$ are positive {\em integers}.  Then the
symplectic form $\omega_{a,b}$ is integral and so the Boothby-Wang
construction \cite{BW} produces a contact form $A_{a,b}$ on the
principal circle bundle $P=P({a,b})$ over $M$ with Chern class $c_1
(P) =[\omega_{a,b}]$.  We will show in Theorem~\ref{main} that for
$a\geq b > 0$ the number of conjugacy classes of maximal tori in the
group of contactomorphisms of $P({a,b})$ is at least the number of
integers $k$ satisfying $0 \leq k < a/b$. Similar results hold for
circle bundles over $\widetilde{\C P}{} ^2$. Since all contact
manifolds $(P(a,b), A_{a,b})$ can be given the structure of a
$K$-contact manifold this answers positively Problem~3 in
\cite{Yamazaki-fbr-j} (in fact the manifolds $P(a,b)$ are
Sasakian-Einstein).

\subsection*{Acknowledgments}
I am grateful to Sue Tolman and Charles Boyer for a number of useful
conversations.  Yael Karshon kindly provided me with a copy of her
preprint. I thank Matthew Ando for answering several of my topology questions.
\section{Preliminaries}

In this section we review contact forms, contact structures, the
Boothby-Wang construction, symplectization, symplectic cones, contact group
actions, contact moment maps, contact and symplectic toric manifolds.

\begin{definition}
Recall that a 1-form $\alpha$ on a manifold $P$ is {\bf contact} if
$\alpha _x \not = 0$ for any $x\in P$, so $\xi = \ker \alpha$ is a
codimension-1 distribution, and if additionally $d\alpha |_\xi$ is
non-degenerate. Thus the vector bundle $\xi \to P$ necessarily has
even-dimensional fibers, and the manifold $M$ is necessarily
odd-dimensional.
A codimension-1 distribution $\xi$ on a manifold $P$ is a
(co-oriented) {\bf contact structure} if it is  a kernel of a
contact form.
\end{definition}
Throughout the paper $\alpha$ will always denote a contact form and
$\xi $ will always denote a contact structure.  We will refer either
to a pair $(P, \alpha)$ or $(P, \xi)$ as a contact manifold.

It is a standard fact that $\xi \subset TP$ is contact iff the
punctured line bundle $\xi^\circ \smallsetminus P$ is a symplectic
submanifold of the cotangent bundle $T^*P$. Here $\xi ^\circ$ denotes
the annihilator of $\xi$ in $T^*P$.  Now suppose $\xi = \ker \alpha$ for
some contact form $\alpha$.  Then $\alpha$, thought of as a map from
$P$ to $T^*P$ is a global nowhere zero section of $\xi^\circ$. Hence
$\xi^\circ \smallsetminus P$ has two components.  Denote one of them
by $\xi^\circ_+$.

\begin{definition}  
A symplectic manifold $(M, \omega)$ is a {\bf symplectic cone} if
there exists a free proper action $\{\rho_t\}$ of the real line on $M$
such that $\rho_t ^*\omega = e^t \omega$ for all $t\in \R$.
The vector field $X$ generating $\{\rho_t\}$ then satisfies $L_X
 \omega = \omega$.  We will refer to $X$ as the {\bf Liouville} vector
 field.  Note that $d\iota (X)\omega = \omega$.
\end{definition}
We recall that if $(M, \omega)$ is a symplectic cone, then the orbit
space $M/\R$ for the action of the real line is naturally a contact
manifold.  We will refer to the contact manifold $M/\R$ as the {\bf
base} of the symplectic cone $(M,\omega)$.

Conversely, if $(P, \xi)$ is a contact manifold then $\xi^\circ_+
\subset T^*P$ is a symplectic cone under the action of $\R$ given by
$\rho_t (p ,\eta) = (p, e^t \eta)$ for all $p\in P$, $\eta \in
(\xi^\circ_+)_p \subset T^*_p P$.  We will refer to $\xi^\circ_+$ as a {\bf
symplectization} of $(P, \xi)$.  If $(M, \omega)$ is a symplectic cone
with the base $(P, \xi)$ then it is not hard to show that $(M,
\omega)$ is symplectomorphic to $\xi^\circ_+$ and that moreover the
symplectomorphism is equivariant with respect to the two $\R$ actions.

If $(M, \omega)$ is a symplectic cone with a Liouville vector field
$X$ and if a Lie group $G$ acts on $M$ preserving $\omega $ and $X$
then $G$ preserves $\iota (X) \omega$.  Hence $\Phi : M \to \fg^*$
given by
\begin{equation} \label{eq-cone-mmap}
\langle \Phi, Y\rangle = \iota (Y_M) \iota (X) \omega
\end{equation}
for all $Y\in \fg$ is a moment map.  Here $Y_M$ denotes a vector
field induced by $Y$ on $M$:  $Y_M (x) = \frac{d}{dt} |_{t=0}
(\exp tY)\cdot x$.  Note that $\Phi (\rho _t (m)) = e^t \Phi
(m)$ for all $t\in \R$, $m\in M$.  

If $(M, \omega)$ is a symplectic manifold and the cohomology class
$[\omega] \in H^2 (M, \R)$ is integral then there exists a principal
circle bundle $\pi: P\to M$ with first Chern class $c_1 (P) =
[\omega]$.  A theorem of Boothby and Wang asserts that moreover there
exists a connection 1-form $A$ on $P$ with $dA = \pi^* \omega$ and
that consequently $A$ is a contact 1-form \cite[Theorem~3]{BW}.  We
will refer to the contact manifold $(P, \xi = \ker A)$ as the {\bf
Boothby-Wang} manifold of $(M, \omega)$.

An action of a Lie group $G$ on a contact manifold $(P, \xi)$ is
contact if it preserves the contact distribution $\xi$.  If
furthermore $G$ is compact and connected then there exists a {\em
$G$-invariant } contact form $\alpha$ with $\ker \alpha = \xi$ (take
any contact form $\alpha_0$ with $\ker \alpha_0 = \xi$ and average it
over $G$).

Suppose now that a Lie group $G$ acts on a manifold $P$ preserving a contact 
form $\alpha$.  We define the corresponding 
{\bf $\alpha$-moment map} $\Psi_\alpha :P
\to \fg^*$   by
\begin{equation}\label{stupid_eq}
\langle \Psi _\alpha (x) , X \rangle = \alpha _x (X_P (x))
\end{equation}
for all $x\in P$ and all $X\in \fg$, where as above $X_P$ denotes the
vector field induced by $X$ on $P$.  One can show that $(P, \alpha,
\Psi_\alpha :P\to \fg^*)$ completely encodes the (infinitesimal) action
of $G$ on $P$.

Note that if $f$ is a $G$-invariant nowhere zero function, then
$\alpha' = f\alpha$ is also a $G$-invariant contact form defining the
same contact structure.  Clearly the corresponding moment map
$\Psi_{\alpha'}$ satisfies $\Psi_{\alpha'} = f\Psi_\alpha$.  Thus the
definition of a contact moment map above is somewhat problematic: it
depends on a choice of an invariant contact form rather then solely on
the contact structure and the action.  Fortunately there is also a
notion of a contact moment map that doesn't have this problem.
Namely, suppose again that a Lie group $G$ acts on a manifold $P$
preserving a contact structure $\xi$ (and its co-orientation).  The
lift of the action of $G$ to the cotangent bundle then preserves a
component $\xi^\circ _+$ of $\xi^\circ \smallsetminus    P$.  The
restriction $\Psi = \Phi|_{\xi_+ ^\circ}$ of the moment map $\Phi$ for
the action of $G$ on $T^*P$ to depends only on the action of the group
and on the contact structure.  Moreover, since $\Phi: T^*P \to
\fg^*$ is given by the formula 
$$
\langle \Phi ( p, \eta ), X\rangle = \langle \eta, X_P (p)\rangle 
$$ 
for all $p\in P$, $\eta \in T^*_p P$ and  $X \in \fg$, we see
that if $\alpha$ is any invariant contact form with $\ker \alpha =
\xi$ then $\langle \alpha^* \Psi ( p,\eta), X\rangle = \langle \alpha ^*
\Phi (p, \eta), X\rangle = \langle \alpha _{p}, X_P (p) \rangle = \langle
\Psi _\alpha (p), X\rangle$.  Here we think of $\alpha$ as a section
of $\xi^\circ_+ \to M$.  Thus $\Psi \circ \alpha = \Psi _\alpha$, that
is, $\Psi = \Phi|_{\xi^\circ_+}$ is a ``universal'' moment map.

\begin{definition}
Let $(P, \xi)$ be a co-oriented contact manifold with an action of a
Lie group $G$ preserving the contact structure $\xi$ and its
co-orientation.  Let $\Psi: \xi^\circ_+ \to \fg^*$ denote the
corresponding moment map.  We define the {\bf moment cone} $C(\Psi)$
to be the set 
$$ 
	C(\Psi) : = \Psi (\xi^\circ_+)\cup \{0\}.  
$$ 
If $\alpha $ is a $G$-invariant contact form with $\xi = \ker
\alpha$ and $\alpha (P) \subset \xi^\circ_+$,  then 
$$
C(\Psi ) = \{t \Psi_\alpha (p) \mid p \in P, \, t \in [0, \infty)\}
= \R^{\geq 0} \Psi_\alpha (P),
$$
where  $\Psi_\alpha :M\to \fg^*$ denote the $\alpha$-moment map.
The moment cone, unlike the $\alpha$-moment map or its image, is an
invariant of the contact group action.
\end{definition}

Recall that a {\bf symplectic toric $G$-manifold} is a triple $(M,
\omega, \Phi:M \to \fg^*)$ where $M$ is manifold with an effective 
Hamiltonian action of a torus $G$ preserving the symplectic form
$\omega$ and satisfying $\dim M = 2\dim G$, and $\Phi$ is a
corresponding moment map.

For this and other reasons we say that an effective contact action of
a torus $G$ on a contact manifold $(P, \xi = \ker \alpha )$ is toric
if the induced action of $G$ on the symplectization $\xi^\circ_+$ of
$P$ makes it a symplectic toric manifold.  Thus a {\bf contact toric
$G$-manifold} is a triple $(P, \xi = \ker \alpha, \Psi_\alpha: P \to
\fg^*)$ where $P$ is a manifold with an effective action of a torus
$G$ on $P$ preserving a contact form $\alpha$ and satisfying $\dim P +
1 = 2 \dim G$, and $\Psi_\alpha : P\to \fg^*$ is the corresponding
$\alpha$-moment map.

\begin{remark}
If a torus $G$ acts effectively and contactly on a contact manifold
$(P, \xi)$ then $\dim G \leq \frac{1}{2} (\dim P + 1)$ (see, for
instance, \cite{LS}).  Thus if $(P, \xi = \ker \alpha, \Psi_\alpha: P
\to \fg^*)$ is a contact toric manifold then $G$ is a maximal torus in
the group of contact diffeomorphisms $\text{ Diff} (P, \xi)$.
\end{remark}

\section{Main theorem}

We are now ready to state the main result of the paper.
\begin{theorem} \label{main}
For every integer $\ell>1$ there exists a set $\{ (P_i, \xi_i =
\ker \alpha_i , \Psi_{\alpha_i} : P_i \to \fh^*)\}_{i=1}^\ell$ of compact
connected contact toric $H$-manifolds which are pairwise
non-isomorphic as contact {\em toric} manifolds but are isomorphic as
contact manifolds.  Hence $(P, \xi) = (P_i, \xi_i)$ has at least $\ell$
non-conjugate (maximal) tori in its group $\text{Diff}(P, \xi)$ of
 diffeomorphisms preserving the contact structure $\xi$.
\end{theorem}

\begin{remark}\begin{enumerate}
\item
We will see that we may take $H$ to be the standard 3-torus $\bbT^3$
and $P$ to be $S^2 \times S^3$.
\item 
By definition $\text{Diff}(P, \xi) = \{ \varphi: P \to P \mid \varphi
\text{ is a diffeomorphism and } d\varphi (\xi) = \xi \}$.
\end{enumerate}
\end{remark}

We now begin our proof of Theorem~\ref{main}.  First let us
spell out exactly what we mean by two contact toric manifolds being
isomorphic.

\begin{definition}
We say that two contact toric $H$-manifolds $(P_i, \xi_i = \ker
\alpha_i, \Psi_{\alpha_i}: P_i \to \fh^*)$, $i=1,2$ are {\em isomorphic} if 
there exists a diffeomorphism $\varphi: P_1 \to P_2$ with $d\varphi
(\xi_1) = \xi _2$ and an isomorphism $\gamma: H\to H$ such that
$$
\varphi (g\cdot p) = \gamma (g) \cdot \varphi (p)
$$
for all $g\in H$, all $p\in P_1$.
\end{definition}
\begin{remark} \label{rmrk3.3'}
If $(P_i, \xi_i = \ker \alpha_i, \Psi_{\alpha_i}: P_i \to \fh^*)$,
$i=1,2$, $\varphi: P_1 \to P_2$ and $\gamma: H\to H$ are as above then
\begin{equation}\label{eq*}
\Psi_{\alpha_2} \circ \varphi = \pm T \circ (e^f \Psi_{\alpha_1})
\end{equation}
for some function $f\in C^\infty (P_1)$.  Here $T= (d\gamma\inv)^*:
\fh^* \to \fh^*$ (the sign in (\ref{eq*}) is $+$ if the contactomorphism 
$\varphi$ preserves the co-orientations of the contact structures).
Hence $\pm T$ maps the moment cone $C(\Psi_1)$ isomorphicly onto the
moment cone $C(\Psi_2)$.  Also, the linear map $T$ preserves the
weight lattice $\Z_H^*$ of $H$, hence $T\in \GL (\Z_H^*)$.  Therefore
in order to prove that two contact toric $H$-manifolds are {\em not}
isomorphic it is enough prove that there does not exist $T\in \GL
(Z_H^*)$ mapping one moment cone onto another.
\end{remark}

Next let me try to explain what makes Theorem~\ref{main} true. As mentioned
in the introduction Yael Karshon constructed for every integer
$\ell>1$ collections $\{(M_i, \omega_i, \Phi_i :M_i \to
\fg^*)\}_{i=1}^\ell$ of symplectic toric $4$-manifolds which are
pairwise non-isomorphic as symplectic {\em toric} manifolds but are
isomorphic as symplectic manifolds \cite{K}.  As we recalled above, by
a theorem of Boothby and Wang for every integral symplectic manifold
$(M, \omega)$ there exists a principal circle bundle $\pi: P\to M$
with a connection 1-form $A$ such that $dA = \pi^* \omega$.  Clearly
if two integral symplectic manifolds are symplectomorphic the
corresponding Boothby-Wang manifolds are contactomorphic.  Therefore
our proof of Theorem~\ref{main} amounts to: choose a Karshon
collection of integral symplectic toric $G$-manifolds $M_i$.  The
corresponding Boothby-Wang manifolds $P_i$ are all contactomorphic.
The actions of $G$ on $M_i$ are covered by actions of extensions of
$G$ by $S^1$ on $P_i$. These extensions are simply $G\times S^1$.
Moreover since $M_i$'s are not isomorphic as symplectic $G$-manifolds,
$P_i$'s should be non-isomorphic as contact $G\times S^1$ manifolds.  We
should be able to ascertain this by looking at their moment cones.
However, I found it simpler to proceed a little differently. Namely
given the moment polytopes of $M_i$'s we will guess what the moment
cones of $P_i$'s should be and then prove that (a) the guess is
correct and (b) there are no elements of $\GL (\Z_G^* \times \Z^*)$
mapping one moment cone into another.

\subsection{Boothby-Wang construction for symplectic toric manifolds}
 
As above we let $\fg^*$ denote the vector space dual of the Lie
algebra of a torus $G$, and $\Z_G^*$ denote the weight lattice of $G$.
A polytope $\Delta \subset \fg^*$ is {\bf Delzant} if for any vertex
$v^*$ of $\Delta$ there exists a basis $\{u_i^*\}$ of $\Z^*_G$ such that
every edge of $\Delta$ coming out of $v^*$ is of the form $\{v^* + t
u_i^* \mid 0\leq t \leq a_i\}$ for some $a_i \geq 0$.  In particular 
$\Delta$ is simple and rational.  Recall

\begin{theorem}[Delzant \cite{D}]\label{thm-Del}  
If $\fg^*$ the dual of the Lie algebra of 
a torus $G$ and $\Delta \subset \fg^*$ is a Delzant polytope then
there exists a unique (compact connected) symplectic toric
$G$-manifold $(M_\Delta, \omega_\Delta, \Phi_\Delta : M_\Delta \to \fg^*)$ 
with $\Phi _\Delta (M_\Delta) = \Delta$.

Conversely if $(M, \omega, \Phi: M\to \fg^*)$ is a (compact connected)
symplectic toric manifold then the moment polytope $ \Phi (M)$ is a
Delzant polytope.
\end{theorem}

We next recall a few well-known facts about (symplectic) toric
manifolds (cf.\ \cite{DavisJan}). Let $(M_\Delta, \omega_\Delta, \Phi:
M_\Delta \to \fg^*)$ be a symplectic toric manifold with moment
polytope $\Delta$.  Then
\begin{enumerate}
\item  The integral cohomology $H^* (M_\Delta, \Z)$ is generated by degree 
2 classes.

\item $H_2 (M_\Delta, \Z)$ is generated by the preimages $\Phi\inv (e)$ 
where $e$'s are the edges of $\Delta$.  Each preimage is a 2-sphere,
and this set of generators is redundant.

\item If $e= \{v^* + tu^* \mid 0\leq t \leq a\}$ is an edge of $\Delta$ and 
$u^* \in \Z_G^*$ is primitive, then $\langle [\omega_\Delta], \Phi\inv
(e)\rangle = a$.

\end{enumerate}

Thus the symplectic form $\omega_\Delta$ is integral iff the edges of
$\Delta$ can be represented by elements of $\Z_G^*$.  Hence if one
vertex of $\Delta$ lies in the weight lattice (and $\omega_\Delta$ is
integral) then all vertices of $\Delta$ lie in the weight lattice.
Therefore it is no loss of generality to assume that if $(M, \omega,
\Phi: M\to \fg^*)$ is a (compact connected) integral symplectic toric 
$G$-manifold then the moment polytope $\Phi (M)$ is an integral
polytope, i.e., that all of its vertices are in the integral lattice $\Z_G^*$.

I showed in \cite{CTM} that an analogue of Theorem~\ref{thm-Del}
holds for symplectic cones.  Let us recall the result.  Let $\fh^*$
denote the dual of the Lie algebra of a torus $H$, $\Z_H^* \subset
\fh^*$ denote its weight lattice and $\Z_H \subset \fh$ denote the
integral lattice.  A cone $C \subset
\fh^*$ is {\bf good} if there exists a (nonempty) set of primitive vectors 
$\mu_1, \ldots , \mu_N \in \Z_H$ such that 
\begin{enumerate}
\item 
$$
C= \bigcap_j \{ \eta \in \fh^* \mid \langle \eta, \mu_j \rangle \geq 0\}
$$
\item The set $\{\mu_1, \ldots, \mu_N\}$ is minimal, i.e.,
$$
C\not = 
\bigcap_{j\not = i} \{ \eta \in \fh^* \mid \langle \eta, \mu_j \rangle \geq 0\}
$$
for any $i$, $0< i \leq N$.

\item Any codimension $k$ face $F$ of $C$ ($k\not = \dim G$) is the 
intersection of exactly $k$ facets, and  the normals to these facets
generate a direct summand of rank $k$ of the integral lattice $\Z_H$.
\end{enumerate}
For example, if $\Delta \subset \fg^*$ is an integral Delzant polytope
then 
$$ 
C_\Delta := 
\{ t (\eta, 1^*) \in \fg^*  \times \R^* \mid \eta \in \Delta, t\geq 0\} 
$$
is a good cone in $\fh^* = \fg^* \times \R^*$.  Here $\{1^*\}$ is a
basis of $\Z^*$ dual to the basis $\{1\}$ of $\Z$, and we think of
$\fh^*$ as the dual of the Lie algebra of $H = G \times \R/\Z$.  We
will refer to $C_\Delta$ as the {\bf standard cone} on the polytope
$\Delta$.

\begin{lemma}
Let $\fh^*$ denote the dual of the Lie algebra of a torus $H$.  Given
a good cone $C\subset \fh^*$ there exists a unique (compact connected)
contact toric manifold $(P_C, \xi_C = \ker \alpha_C, \Psi_{\alpha_C}:
P_C \to \fh^*)$ with moment cone $C$.

Equivalently there exists a unique symplectic toric cone (with a
compact connected base) $(M_C, \omega_C, \Phi_C :M_C \to \fh^*)$ such
that $\Phi _C (M_C) = C \smallsetminus \{0\}$.  Additionally the
fibers of $\Phi_C$ are connected.
\end{lemma}
\begin{proof}
See Theorem~2.8 in \cite{CTM}.
\end{proof}

\begin{lemma}\label{BW}
Let $\fg^*$ denote the dual of the Lie algebra of a torus $G$,
$\Delta\subset \fg^*$ be an integral Delzant polytope, $C= C_\Delta
\subset \fg^* \times \R^*$ the standard cone on $\Delta$, $(M_\Delta,
\omega_\Delta, \Phi_\Delta: M_\Delta \to \fg^*)$ the compact connected
integral symplectic toric manifold with moment polytope $\Delta$ and
$(M_C, \omega_C, \Phi_C :M_C \to \fg^*\times \R^* )$ the symplectic
toric cone corresponding to the cone $C$.  Let $X$ denote the
Liouville vector field on the symplectic cone $(M_C, \omega_C)$.  Then
\begin{enumerate}
\item $S^1 = \{1\} \times \R/\Z \subset G \times \R/\Z$
acts freely on $P_C: = \Phi_C\inv (\fg^* \times \{1^*\})$ making it a
principal $S^1$ bundle over $P_C/S^1$.

\item The quotient $P_C/S^1$ is naturally a symplectic toric $G$-manifold 
(isomorphic to) $(M_\Delta, \omega_\Delta, \Phi_\Delta)$.

\item $\alpha_C := (\iota(X) \omega_C)|_{P_C}$ is a connection 1-form on the 
principal $S^1$-bundle $\pi_C :P_C \to M_\Delta$. 

\item ${\pi_C}^* \omega_\Delta = d\alpha_C$.
\end{enumerate}
Hence $(P_C, \alpha_C)$ is the Boothby-Wang manifold of $(M_\Delta,
\omega_\Delta)$, and $(M_C,  \omega_C)$ is the symplectization of 
$(P_C, \alpha_C)$.
\end{lemma}

\begin{proof}
The action of $S^1$ at a point $x\in M_C$ is free iff $S^1$ intersects
trivially the isotropy group $H_x$ of $x$ in $H= G\times \R/\Z$.  Since
$M_C$ is $H$-toric, $H_x$ is connected hence a torus.  Now two subtori
$T_1$, $T_2$ of a torus $K$ intersect trivially iff their integral
lattices $\Z_{T_1}$, $\Z_{T_2}$ intersect trivially and $\Z_{T_1} +
\Z_{T_2}= \Z_{T_1} \oplus \Z_{T_2}$ is a direct summand of the integral 
lattice $\Z_K$ of $K$.  The Lie algebra $\fh_x$ of $H_x$ is the
annihilator of the face $F$ of $C$ containing $\Phi_C (x)$ in its
interior.  The integral lattice of $H_x$ is generated by the normals
to the facets of $C$ intersecting in $F$.

Let us consider the worst case scenario: $\dim H_x = \dim H -1$,
i.e., the face $F$ is the ray $\R^{>0} (v^*, 1^*)$ through a vertex
$(v^*, 1^*)$ of $\Delta\times \{1^* \}$.  The general case then
follows easily from this special case.  If $\{u^*_i\}$ is a basis of
$\Z_G^*$ spanning the edges of $\Delta$ coming out of the vertex
$v^*$, then the dual basis $\{u_i\}\subset \Z_G$ consists of the
normals to the facets of $\Delta$ meeting at $v^*$.  The vectors
$\{(u_i, - \langle v^*, u_i\rangle 1)\}\subset \Z_G \times \Z = \Z_H$
are then the normal vectors to the facets of $C$ that intersect in the
ray $\R^{>0} (v^*, 1^*)$.  Clearly the set $\{(u_i, - \langle v^*,
u_i\rangle 1)\} \cup \{ (0, 1)\}$ is a basis of $\Z_G \times Z$.  We
conclude that the action of the $S^1 $ on $M_C$ is free.

Since $S^1$ acts freely, 1 is a regular value of $f:= \langle \Phi_C,
(0, 1)\rangle$.  Hence $P_C = f\inv (1) = \Phi_C \inv (\fg^* \times
\{1^*\}) = \Phi_C \inv ((\fg^* \times \{1^*\}) \cap C ) = \Phi_C \inv
(\Delta \times \{1^*\})$ is a principal $S^1$ bundle over
$P_C/S^1$. It is connected since the fibers of $\Phi_C$ are connected.
By the symplectic reduction theorem of Marsden - Weinstein and  Meyer,
$P_C/S^1$ is naturally a symplectic manifold.  Moreover, since the
actions of $S^1$ and $G$ commute, the action of $G$ on $M_C$ induces a
Hamiltonian action of $G$ on $P_C/S^1$.  Furthermore the restriction
of $\Phi_C$ to $P_C$ descends to a corresponding $G$ moment map $\Phi$
on $P_C/S^1$ (provided we identify $\fg^*$ and $\fg^* \times
\{1^*\}$).  Clearly $\Phi (P_C/S^1) = \Delta \times \{1^*\}$.
Therefore by the uniqueness part of Delzant's classification of
compact symplectic toric manifolds (Theorem~\ref{thm-Del} above)
$P_C/S^1$ is $(M_\Delta, \omega_\Delta, \Phi_\Delta)$.

Finally we argue that $\omega_\Delta$ is the curvature of a connection
1-form on the principal $S^1$-bundle $\pi_C: P_C \to M_\Delta$.  By
construction ${\pi_C}^* \omega_\Delta = \omega_C |_{P_C}$.  Let $Y$
denote the vector field on $M_C$ generating the action of the $S^1$.
Then $\langle \Phi_C, (0, 1)\rangle = \iota (Y) \iota (X)\omega_C$
(cf.\ equation (\ref{eq-cone-mmap})).  Hence $1 = \iota (Y)\alpha_C$
where $\alpha _C : = \iota (X)\omega_C)|_{P_C}$.  It follows that
$\alpha_C$ is a connection on $\pi_C : P_C \to M_\Delta$.  Moreover,
since $d\iota(X) \omega_C = \omega _C$, we have $d\alpha_C = (d
\iota(X)\omega_C)|_{P_C} = \omega_C |_{P_C} = {\pi_C}^* \omega_\Delta$.
\end{proof}

Next we apply Lemma~\ref{BW} to construct certain contact toric
5-manifolds out of symplectic toric Hirzebruch surfaces.  It will be
convenient to take the standard $n$-torus $\bbT^n$ to be $\R^n /\Z^n$
and to identify the Lie algebra of $\bbT^n$ with $\R^n$.  By using the
standard basis of $\R^n$ we identify the dual of the Lie algebra of
$\bbT^n$ with $\R^n$ and the weight lattice of $\bbT^n$ with $\Z^n$.

\begin{definition}
Following Karshon \cite{K} we define the {\bf standard Hirzebruch
trapezoid} $\Delta(a,b, m) $ to be the quadrilateral in $\R^2$ with
vertices $(0,0)$, $(0,b)$, $(a - \frac{m}{2}b, b)$ and $(a+
\frac{m}{2}b, 0)$ where $m$ is a non-negative integer and $b> 0$, $a >
\frac{m}{2}b$ are real numbers.  
\end{definition}
By Theorem~\ref{thm-Del} there exists a symplectic 4-manifold
$(M(a,b,m), \omega_{a,b,m})$ with an effective Hamiltonian action of
$\bbT^2$ such that $\Delta (a,b, m)$ is the image of $M(a, b, m)$
under the corresponding moment map.  One can show that the manifold
$M(a,b,m)$ is a Hirzebruch surface. In particular it is diffeomorphic
to either $S^2 \times S^2$ or to $\widetilde{\C P} {}^2$, depending on
the values of $a$, $b$, and $m$.
\begin{definition}
We define the {\bf standard Hirzebruch cone} $C(a,b,m)$ to be the
standard cone on the Hirzebruch trapezoid $\Delta(a,b,m)$: 
$$ 
C(a,b, m) = \{ t (x_1, x_2, 1) \in \R^3 \mid t \geq 0, (x_1, x_2 ) \in \Delta
(a,b, m)\}.  
$$
\end{definition}
Now suppose that $b$ and $ a - \frac{m}{2}b$ are integers.  Then the
Hirzebruch trapezoid $\Delta(a, b, m)$ is integral. It follows from
Lemma~\ref{BW} that $C(a,b,m)$ is the moment cone of the Boothby-Wang
manifold $(P(a,b,m),\ker A_{a,b,m})$ of the integral symplectic
manifold $(M(a,b,m), \omega _{a,b,m})$.

\begin{proposition}\label{prop3.10}
Let $C(a,b,m)$ and $C(a',b',m')$ be two standard Hirzebruch cones with
$b,b', a - \frac{m}{2}b, a' - \frac{m'}{2}b' \in \Z$.  If there is
$T\in \GL(3, \Z)$ with $T(C(a,b, m)) = C(a',b',m')$ then either
$(a,b,m) = (a', b',m')$ or $m= m' = 0$ and $a=b'$, $b=a'$.
\end{proposition}

\begin{proof}
A Hirzebruch trapezoid $\Delta (a,b,m)$ has the following property:
if $v_0$ is a vertex of $\Delta (a,b,m)$ and $v_1$, $v_2$ are two
adjacent vertices then there is a basis $\{u_1, u_2\}$ of $\Z^2$ such
that $v_1 -v_0 = t_1 u_1$, $v_2 -v_0 = t_2 u_2$ for some $t_1, t_2
>0$.  Depending on which three vertices we picked the set $\{t_1,
t_2\}$ is either $\{b, a - \frac{m}{2} b\}$ or $\{ b, a + \frac{m}{2}
b\}$.

If $w\in \Z^3$ is a primitive vector then for any $T\in \GL (3, \Z)$
the vector $Tw$ is also primitive.  Now suppose $T (C(a,b, m))=
C(a',b',m')$.  Then $T$ maps the edge of $C(a,b,m)$ through $(0,0, 1)$
to an edge of $C(a', b', m')$ say through a vertex $v_0$ of
$\Delta(a', b' , m') \times \{1\}\subset \R^3$.  Since both $T(0,0,1)$
and $v_0$ are primitive, we have $T(0,0,1) = v_0$.  Similarly $T$
maps the vectors $( 0, b, 1)$ and $(a+ \frac{m}{2} b, 0, 1)$ to
vertices $v_2$, $v_1$ of $\Delta(a', b' , m') \times \{1\}$ adjacent
to $v_0$.  It follows from the remark at the start of the proof that
there are vectors $u_1, u_2 \in \Z^3$ such that $\{v_0, u_1, u_2\}$ 
is a basis of $\Z^3$ and 
$v_1 - v_0 = t_1 u_1$, $v_2 - v_1 = t_2 u_2$ where 
\begin{equation} \label{eq1}
\{t_1, t_2\} = 
\{b', a' - \frac{m'}{2} b'\} \text{ or }\,\{ b', a' + \frac{m'}{2} b'\}.
\end{equation}
Now $T( 0, 1,  0) = T (\frac{1}{b} (( 0, b, 1) - (0, 0, 1)))
= \frac{1}{b} (v_2 - v_0)  = \frac{t_2}{b} u_2$.  Similarly $T(1, 0, 0)= 
\frac{t_2}{a + \frac{m}{2}b} u_1$.  Since both $(0,1, 0)$ and $u_2$ are 
primitive in $\Z^3$, $\frac{t_2}{b} = \pm 1$.  Since $t_2 > 0$ we get
\begin{equation}\label{eq3}
t_2 = b.
\end{equation}
By the same argument 
\begin{equation} \label{eq4}
t_1 = a + \frac{m}{2} b.
\end{equation}
Finally let $v_3$ denote the remaining vertex of
$\Delta (a',b',m')\times \{1\}$.  Then $v_3 = T (a-\frac{m}{2} b, b, 1)
= (a -\frac{m}{2}b) T(1,0,0) + b T (0,1,0) + T (0,0,1)$.  Hence
\begin{equation} \label{eq2}
v_3=  (a -\frac{m}{2}b) u_1 + b u_2 + v_0.
\end{equation}
Equation (\ref{eq1}) gives us four cases to compare (\ref{eq2}) with.
For example suppose $t_1 = a' + \frac{m'}{2} b'$ and $t_2 = b'$.  Then
$$
 v_3 = v_0 + (a' - \frac{m'}{2}b') u_1 + b' u_2.  
$$ 
Comparing the above equation with equation (\ref{eq2}) and using
equations (\ref{eq3}) and (\ref{eq4}) we get:
\begin{align*}
a+ \frac{m}{2} b & = a' + \frac{m'}{2} b'\\
a- \frac{m}{2} b & = a' - \frac{m'}{2} b'\\
b& = b' .
\end{align*}
Hence $a= a'$, $b= b'$, $m=m'$.  The proposition follows by examining the
remaining three cases.
\end{proof}
Combing Proposition~\ref{prop3.10} with Remark~\ref{rmrk3.3'} we see
that the contact toric $\bbT^3$ manifolds $(P(a,b, m), \ker
A_{a,b,m})$ are isomorphic as contact {\em toric} manifolds if and
only if either $(a,b,m) = (a', b',m')$ or ($m= m' = 0$ and $a=b'$,
$b=a'$).  On the other hand, by Lemma~3 of \cite{K} the manifold
$(M(a,b,m),
\omega_{a,b,m})$ is symplectomorphic to $(M(a, b, m+2), \omega_{a, b,
m+2})$. Hence the 5-manifolds $(P(a, b, m), \ker A_{a,b,m})$ and
$(P(a, b, m+2), \ker A_{a, b,m+2})$ are contactomorphic.  Now $P(a,
b,m)$ is well-defined as long as $\frac{a}{b} > \frac{m}{2} \geq 0$.
Let $k$ be the largest positive integer with $\frac{a}{b} >
\frac{k}{2}$.  Then $P(a,b, k)$, $P(a, b, k-2)$, \ldots are contact
toric $\bbT^3$-manifolds which are isomorphic as contact manifolds but
not as contact {\em toric} manifolds.  Note that there are $\ell$ of these
manifolds where $k = 2\ell-2$ or $2\ell -1$.  This finishes our proof of
Theorem~\ref{main}.

\begin{remark}
By a theorem of Hatakeyama \cite{H} the Boothby-Wang manifold of an
integral Kaehler manifold is Sasakian.  Hence the manifolds $P(a,b,
m)$ considered above are all Sasakian.  In particular they are all
$K$-contact.   Therefore Theorem~\ref{main} gives a positive answer 
to Problem~3 in \cite{Yamazaki-fbr-j}.
\end{remark}

\section{Contact structures on $S^2 \times S^3$ }

In this section we examine the manifolds constructed in the course of
the proof of Theorem~\ref{main} in more details.  We identify
precisely what some of these manifolds are: for $a$ and $b$ relatively
prime the manifolds $P(a, b, 0)$ turn out to be $S^2 \times S^3$.  We
also identify the contact structures we have constructed as complex
vector bundles.  We do not know  whether these contact
structures which happen to be isomorphic as vector bundles are in fact
isomorphic as contact structures.

\begin{proposition}\label{thm2}
Suppose $a, b$ are two relatively prime integers with $a> b \geq 1$.
Let $(P(a,b, 0), \xi_{a, b} = \ker A_{a,b, 0}) $ denote the
corresponding contact toric 5-manifold defined above.  It is the
Boothby-Wang manifold of $(\C P^1 \times \C P^1, \omega_{a, b})$.
\begin{enumerate}
\item  $P(a, b, 0)$ is diffeomorphic to $S^2 \times S^3$.  
\item As a complex vector bundle $\xi_{a, b}$ is uniquely determined by the 
difference $a-b$.
\end{enumerate}
\end{proposition}

\begin{proof} 
The first part of the proposition is an observation of Wang and Ziller
\cite[(2.3)]{WZ}.  Since the manifold $P(a, b, 0)$ is the Boothby-Wang
manifold of $(\C P^1 \times \C P^1, \omega_{a,b})$ the orbit map $\pi:
P(a,b, 0) \to \C P^1 \times \C P^1$ is a principal circle bundle with
Chern class $[\omega_{a,b}]$.  Principal torus bundles over products
of Kaehler-Einstein manifolds were studied extensively by Wang and
Ziller, who showed that these bundles carry Einstein metrics (in the
notation of \cite{WZ} $P(a,b, 0)$ is $M^{1,1}_{a,b}$).  In particular
they showed that if $a$ and $b$ are relatively prime then $P(a,b, 0)$
is diffeomorphic to $S^2 \times S^3$ (cf.\
\cite[(2.3)]{WZ}).\footnote{The rough idea of the proof is to first
show that $P(a,b, 0)$ is simply connected and spin.  It is easy to
compute that $H^2 (P(a,b, 0), \Z)$ is $\Z$.  Then Smale's
classification of 5-dimensional simply connected spin manifolds implies
that $P(a,b,0) = S^2 \times S^3$.}  

We now prove part 2 of the theorem. Suppose $(a', b')$ is another pair
of relatively prime integers with $a' - b' = a-b$.  We claim that
$\xi_{a,b}$ and $\xi_{a',b'}$ are isomorphic as complex vector
bundles.  In the course of the proof of \cite[(2.3)]{WZ} Wang and
Ziller showed that if one denotes the generator of $H^2 (P(a,b, 0),
\Z)$ by $z$ and the obvious generators of $H^2 (\C P^1 \times \{*\},
\Z)$ and $H^2 (\{*\}\times \C P^1, \Z)$ by $x$ and $y$ respectively,
then $\pi^* x = az$ and $\pi^* y = -bz$.  Consequently the first Chern
class of the line bundle $L_1^{a,b} := \pi^* T (\C P^1\times \{*\})$
is $\pi^* (2x) = 2az$, and similarly the first Chern class of $L_2
^{a, b} := \pi^* T (\{*\}\times \C P^1)$ is $-2bz$.  Since the contact
distribution $\xi_{a,b}$ is a connection on the bundle $\pi: P(a,b, 0)
\to \C P^1 \times \C P^1$, we have $\xi_{a,b} = \pi^* T (\C P^1 \times
\C P^1)$.  Hence $\xi_{a,b}= L_1^{a,b} \oplus L_2 ^{a, b}$.  The
distribution $\xi_{a',b'}$ is also a direct sum of line bundles:
$\xi_{a', b'} = L_1^{a',b'} \oplus L_2 ^{a', b'}$; additionally  $c_1
(L_1^{a',b'}) = 2a'z$ while $c_1 (L_2 ^{a', b'}) = -2b'z$.  The
projection $p : S^2 \times S^3 \to S^2$ induces an isomorphism on the
level of second cohomology.  Therefore all the line bundles over
$P(a,b, 0) = S^2 \times S^3 = P(a', b', 0)$ are pull-backs by $p$ of
line bundles over $S^2$, and consequently $\xi_{a, b}$ and $\xi_{a',
b'}$ are pull-backs by $p$ of rank 2 vector bundles over $S^2$.  Two
complex rank $n$ vector bundles over $S^2$ are isomorphic if and only
if their determinant bundles are isomorphic.  It follows that $\xi_{a,
b} $ is isomorphic to $\xi_{a', b'}$ iff $a-b = a' -b'$.
\end{proof}

We end the paper with two question.
\begin{question}
Suppose $a, b$ and $a', b'$ are two different pairs of relatively
prime integers with $a > b \geq 1$, $a' > b' \geq 1$ and $a-b = a'
-b'$.  Are the two contact manifolds $(P(a, b, 0), \xi_{a,b})$ and
$P(a', b', 0), \xi_{a', b'})$ contactomorphic?
\end{question}

\begin{question}
What does Theorem~\ref{main} tell us about the topology of the group
of contactomorphisms $\text{Diff}(P(a,b,m), \ker A_{a,b,m})$?  For
example, do the non-conjugate tori define different homology classes? 
\end{question}


\begin{thebibliography}{WWWW}

\bibitem[Ab]{Ab} M. Abreu, 
Topology of symplectomorphism groups of $S\sp 2\times S\sp 2$,
{\em Invent.\ Math.} {\bf 131} (1998), no. 1, 1--23. 

\bibitem[AbMc]{AbMc} M. Abreu and D. McDuff, Topology of 
symplectomorphism groups of rational ruled surfaces, {\em J.\ Amer.\
Math.\ Soc.} {\bf 13} (2000), no. 4, 971--1009

\bibitem[BW]{BW} W.M. Boothby and H.C. Wang, On contact manifolds,
{\em Ann.\ of Math.\ (2)} {\bf 68} (1958), 721--734. 

\bibitem[DJ]{DavisJan} M.W. Davis and T. Januszkiewicz, Convex polytopes, 
Coxeter orbifolds and torus actions, {\em Duke Math.\ J.} {\bf 62}
(1991), no. 2, 417--451.

\bibitem[D]{D} T. Delzant,  Hamiltoniens p\'eriodiques et images convexes de 
l'application moment, {\em Bull.\ Soc.\ Math.\ France} {\bf 116}
(1988), no. 3, 315--339.

\bibitem[G]{Gr} M. Gromov, Pseudoholomorphic curves in symplectic manifolds, 
{\em Invent.\ Math.}  {\bf 82}, (1985), 307--347.

\bibitem[H]{H} Y. Hatakeyama,  
Some notes on differentiable manifolds with almost contact structures,
{\em Osaka Math.\ J.\ (2)} {\bf 15} 1963, 176--181.

\bibitem[K]{K} Y. Karshon, Maximal tori in the symplectomorphism groups of 
Hirzebruch surfaces, preprint, {\tt
http://www.ma.huji.ac.il/$\sim$karshon/papers/}.

%\bibitem[L1]{pi1} E. Lerman, Homotopy groups of $K$-contact toric manifolds.

\bibitem[L2]{CTM} E. Lerman, Contact toric manifolds, 
{\tt http://xxx.lanl.gov/abs/math.SG/0107201}.

\bibitem[L3]{L-IsJM} E. Lerman, Contact cuts, {\em Israel J.\ Math}, 
{\bf 124} (2001), 77--92;
{\tt http://xxx.lanl.gov/abs/math.SG/0002041}


\bibitem[LS]{LS} E. Lerman and N. Shirokova, Completely integrable torus 
actions on symplectic cones, {\em Math.\ Res.\ Lett.} {\bf 9} (2002).

\bibitem[WZ]{WZ} M.Y. Wang and W. Ziller, Einstein metrics on principal torus 
bundles, {\em J.\ Differential Geom.} {\bf 31} (1990), no. 1,
215--248.

\bibitem[Y1]{Yamazaki-fbr-j} T. Yamazaki, A construction of $K$-contact 
manifolds by fiber join, {\em Tohoku Math.\ J.} {\bf 51} (1999),
433--446.


\end{thebibliography}
\end{document}